\documentclass[a4paper,12pt,reqno]{amsart}
\usepackage{amsmath,a4,amsfonts,amssymb,amsthm,enumerate,fleqn,breqn, dsfont,,stackrel,mathtools}
\usepackage[a4paper,left=1in, right=1in, top=1.4in, bottom=1.1in]{geometry}
\usepackage{tikz}

\usepackage{enumerate}

\usepackage{enumitem}

\def\RR{\mathbb{R}}
\def\NN{\mathbb{N}}
\def\HH{\mathfrak{Q}}

\def\cc#1{\{#1\}}
\def\pp#1{\|#1\|}
\def\ra{\rangle}
\def\la{\langle}
\def\nh{\mathcal{H}}

\def\E{\mathcal{E}}
\def\H{V_R(\HH)}
\def\E{\ell_2(\HH)}

\theoremstyle{plain}
\newtheorem{theorem}{\bf Theorem}[section]

\theoremstyle{remark}
\newtheorem{definition}[theorem]{\bf Definition}

\newtheorem{example}[theorem]{\bf Example}
\newtheorem{remark}[theorem]{\bf Remark}

\newtheorem{corollary}[theorem]{\bf Corollary}

\title[F\lowercase{rames in Quaternionic} H\lowercase{ilbert Space}]{Frames in Quaternionic Hilbert Spaces}
\author[S\lowercase{harma and}  G\lowercase{oel} ]{S.K. S\lowercase{harma and} S\lowercase{hashank} G\lowercase{oel}  }

\address[S.K. Sharma]{D\lowercase{epartment} \lowercase{of} M\lowercase{athematics},
	K\lowercase{irori } M\lowercase{al} C\lowercase{ollege},
	U\lowercase{niversity of} D\lowercase{elhi}, D\lowercase{elhi~110~007}, INDIA.}
\email{\lowercase{sumitkumarsharma@gmail.com}}

\address[Shashank Goel]{D\lowercase{epartment} \lowercase{of} M\lowercase{athematics},
	A\lowercase{mity} I\lowercase{nstitute of} A\lowercase{pplied}  S\lowercase{ciences},
	A\lowercase{mity} U\lowercase{niversity}, N\lowercase{oida}, U.P-201301, INDIA.}
\email{\lowercase{goel.shashank25@gmail.com}}

\subjclass[2010]{42C15, 42A38} \keywords{Frame, Quaternionic Hilbert spaces.} 

\begin{document}
\maketitle \baselineskip12pt

\begin{abstract}
In this paper, we introduce and study the frames in  separable quaternionic Hilbert spaces. Results on the existence of frames in quaternionic Hilbert spaces have been given. Also, a characterization
of frame in quaternionic Hilbert spaces in terms of frame operator is given. Finally, a
Paley-Wiener type perturbation result for  frames in quaternionic Hilbert space has been obtained.
\end{abstract}

\section{Introduction}

Duffin and Schaeffer \cite{DS}   indroduced \emph{frames for Hilbert spaces} while working on some deep problems in
non-harmonic Fourier series. They gave the following definition:

``Let $\nh$ be a Hilbert space. Then $\{x_n\}_{n\in\NN}\subset \nh$ is said to be a \emph{frame}
for  $\nh$ if there exist finite constants $A$ and $B$ with
 $0<A\le B$ such that
\begin{eqnarray}\label{1.1}
A\|x\|^2\le \sum\limits_{n=1}^\infty |\langle
x,x_n\rangle|^2\le B\|x\|^2, \ \ \text{for all} \ x\in \nh."
\end{eqnarray}
The positive constants $A$ and $B$, respectively, are called lower
and upper frame bounds for the frame $\cc{x_n}_{n\in\NN}$.
The inequality (\ref{1.1}) is called the \emph{frame inequality} for the frame
$\{x_n\}_{n\in\NN}$.
A frame $\cc{x_n}_{n\in\NN}$ in $\nh$ is said to be
\begin{itemize}[leftmargin=.35in]
	\item \emph{{tight}} if it is possible to choose   $A=B$.
	\item \emph{Parseval} if it is a tight frame with $A=B=1$.
\end{itemize}

Frame theory began to spread among researchers when in 1986,
Daubechies, Grossmann
and Meyer published a fundamental paper \cite{DGM} and observed that frames can be also used to provide the series
expansions of functions in $L^2({\RR})$.
The main property of frames which makes them so useful is their redundancy, due to which,  representation
of signals using frames is advantageous over basis expansions in a variety of practical applications.

Frame works as an important tool in the study of  signal and image processing \cite{BCE}, filter bank theory \cite{BHF},
wireless communications \cite{HP} and sigma-delta quantization \cite{BPY}. For more literature on
frame theory, one may refer to \cite{C1, CH2, CK}.

Recently, Khokulan, Thirulogasanthar and Srisatkunarajah \cite{KTS} introduced and studied frames for finite dimensional quaternionic Hilbert spaces. Sharma and Virender \cite{SV} study some different types of dual frames of a given frame in a  quaternionic Hilbert space and gave various types of reconstructions with the help of  dual frame.  In this paper, we will  introduce and study the frames in separable quaternionic Hilbert spaces . Results on the existence of frames in quaternionic Hilbert spaces have been given. Also, a characterization
of frame in quaternionic Hilbert spaces in terms of frame operator is given. Finally, a
Paley-Wiener type perturbation result for  frames in quaternionic Hilbert space has been obtained.
\section{Quaternionic Hilbert space}
\setcounter{equation}{0}
\fontsize{12}{14}

As the quaternions are non-commutative in nature therefore there are two different types of quaternionic Hilbert spaces exist, the left  quaternionic Hilbert space and the right quaternionic Hilbert space depending on positions of quaternions. In this section, we will study some basic notations about the algebra of quaternions, right quaternionic Hilbert space and operators on right quaternionic Hilbert spaces. 

Throughout this paper, we will denote $\mathfrak{Q}$ to be a non-commutative field of quaternions,  $I$ be a non empty set of indicies, $\H$ be a separable right quaternionic Hilbert space,  by the	term ``right linear operator", we mean a ``right $\HH$-linear operator" and $\mathfrak{B}(\H)$ denotes the set of all bounded (right $\HH$-linear) operators of $\H$:
\begin{eqnarray*}
	\mathfrak{B}(\H) := \{T : \H\rightarrow\H: \|T\|<\infty\}.
\end{eqnarray*}

\bigskip

The non-commutative field of quaternions $\mathfrak{Q}$ is a four dimensional real algebra with unity. In $\HH$, $0$ denotes the null element and $1$ denotes the identity with respect to multiplication. It also includes three so-called imaginary units, denoted by $i,j,k$. i.e.,
\begin{eqnarray*}
	\mathfrak{Q}=\cc{x_0+x_1i +x_2j +x_3k \ :\ x_0,\ x_1,\ x_2,\  x_3\in \RR}
\end{eqnarray*}
where $i^2=j^2=k^2=-1; \ ij=-ji=k; \ jk=-kj=i$ and $ki=-ik=j$. For each quaternion $q=x_0+x_1i +x_2j +x_3k \in \mathfrak{Q}$, define conjugate of $q$
denoted by $\overline{q}$ as
\begin{eqnarray*}
	\overline{q}=x_0-x_1i -x_2j -x_3k \in \mathfrak{Q}.
\end{eqnarray*}
If $q=x_0+x_1i +x_2j +x_3k$ is a quaternion, then $x_0$ is called the real part of $q$ and $x_1i +x_2j +x_3k$ is called the imaginary part  of $q$. The modulus of a quaternion $q=x_0+x_1i +x_2j +x_3k$ is defined as
\begin{eqnarray*}
	|q|=(\overline{q}q)^{1/2} = (q\overline{q})^{1/2}= \sqrt{x_0^2 +x_1^2 +x_2^2 +x_3^2 }.
\end{eqnarray*}
For every non-zero quaternion $q=x_0+x_1i +x_2j +x_3k \in \mathfrak{Q}$, there exists a unique inverse $q^{-1}$ in $\mathfrak{Q}$ as
\begin{eqnarray*}
	q^{-1}=\dfrac{\overline{q}}{|q|^2 } = \dfrac{x_0-x_1i -x_2j -x_3k }{\sqrt{x_0^2 +x_1^2 +x_2^2 +x_3^2 }}.
\end{eqnarray*}

\begin{definition}
	A \textit{right quaternionic vector space} $\mathds{V}_R(\HH)$ is a linear vector space under right scalar multiplication over the field of quaternionic $\HH$, i.e.,
	\begin{eqnarray}\label{2.1}
	\mathds{V}_R(\HH)\times\HH &\rightarrow& \mathds{V}_R(\HH) \nonumber\\ 
	(u,q)&\rightarrow& uq
	\end{eqnarray}
	and for each $u, v\in\mathds{V}_R(\HH)$ and $p, q\in\HH$, the right scalar multiplication (\ref{2.1}) satisfying the following properties: 
	\begin{eqnarray*}
		&&(u+v)q=uq+vq\\
		&&u(p+q)=up+uq\\
		&&v(pq)=(vp)q.
	\end{eqnarray*}
\end{definition}

\begin{definition}\label{qhs}
	A \textit{right quaternionic pre-Hilbert space} or \textit{right quaternionic inner product space} $\mathds{V}_R(\HH)$ is a right quaternionic vector space together with the binary mapping
	$\langle . | . \rangle : \mathds{V}_R(\HH) \times \mathds{V}_R(\HH) \to \mathfrak{Q}$ (called the \textit{Hermitian quaternionic inner product})
	which satisfies following properties:
	\begin{enumerate}[label=(\alph*)]
		\item $\overline{\langle v_1 | v_2 \rangle} = \langle v_2 | v_1 \rangle$  for all $v_1, v_2 \in \mathds{V}_R(\HH)$.
		\item $\langle v | v \rangle > 0 \ \text{if} \ \ v\ne 0$.
		\item $\langle v | v_1 + v_2 \rangle = \langle v | v_1 \rangle + \langle v | v_2 \rangle $  for all $v, v_1, v_2 \in \mathds{V}_R(\HH)$
		\item $\langle v | uq \rangle = \langle v | u \rangle q $  for all $v, u \in \mathds{V}_R(\HH)$ and $ q \in \mathfrak{Q}$.
	\end{enumerate}
\end{definition}

In view of Definition \ref{qhs}, a right quaternionic inner product  space $\H$  also have  the following properties:
\begin{enumerate}[label=(\roman*)]
	\item $\langle qv | u \rangle = \overline{q}\langle v | u \rangle $  for all $v, u \in \mathds{V}_R(\HH)$ and $ q \in \mathfrak{Q}$
	\item $pv_1 +qv_2\in \mathds{V}_R(\HH)$, for all $v_1, v_2 \in \mathds{V}_R(\HH)$ and $p, q \in \mathfrak{Q}$
\end{enumerate}

Let $\mathds{V}_R(\HH)$ be right quaternionic inner product space with the Hermitian inner product $\langle .|.\rangle$. ~Define the quaternionic norm $\|.\|:\mathds{V}_R(\HH)\rightarrow\RR^+$ on $\mathds{V}_R(\HH)$ by
\begin{eqnarray}\label{2.2}
\|u\|=\sqrt{\langle u|u\rangle},\ u\in\mathds{V}_R(\HH).
\end{eqnarray}

\begin{definition}
	The right quaternionic pre-Hilbert space  is called a \textit{right quaternionic Hilbert space}, if it is complete with respect to the norm (2.2) and is denoted by $\H$.
\end{definition}

\begin{theorem}[The Cauchy-Schwarz Inequality]\cite{GMP} 
If $\H$ is a right quaternionic Hilbert space then  
\begin{eqnarray*}
|\langle u|v\rangle|^2\leq\langle u|u\rangle\langle v|v\rangle, \ \ \text{for all}\ \ u, v \in \H.
\end{eqnarray*}
Moreover, a norm defined in (\ref{2.2}) satisfy the following properties:
\begin{enumerate}[label=(\alph*)]
	\item $\|uq\|=\|u\||q|$, for all $u\in\H$ and $q\in\HH$.
	\item $\|u+v\|\leq\|u\|+\|v\|$, for all $u, v\in\H$.
	\item $\|u\|=0$ for some $u\in\H$, then $u=0$.
\end{enumerate}
\end{theorem}

\noindent
	For the non-commutative field of quaternions $\HH$,  define the space $\ell_2(\HH)$ by
	\begin{eqnarray*}
		\ell_2(\HH) = \bigg\{\{q_i\}_{i\in I}\subset \HH : \ \sum_{i\in I}\|q_i\|^2<+\infty \bigg\}
	\end{eqnarray*}
	under right multiplication by quaternionic scalars together with   the  quaternionic inner product on $\ell_2(\HH)$ defined as
	\begin{eqnarray}\label{2.3}
	\langle p|q\rangle=\sum_{i\in I}\overline{p_i}q_i,\ p=\{p_i\}_{i\in I}\ \text{and}\ q=\{q_i\}\in\ell_2(\HH).
	\end{eqnarray}
It is easy to observe that $\ell_2(\HH)$ is a right quaternionic Hilbert space with respect to quaternionic inner product (\ref{2.3}).

\begin{definition}[\cite{GMP}] Let $\H$ be a right quaternionic Hilbert Space and  $S$ be a subset of $\H$. Then, define the set:
	\begin{itemize}[leftmargin=.35in]
	\item	$S^{\bot}=\{v\in\H:\langle v|u\rangle=0\ \forall\ u\in S\}.$
	\item $\langle S\rangle$ be the right $\HH$-linear subspace of $\H$ consisting of all finite right $\HH$-linear combinations of elements of $S$.
	\end{itemize}
 \end{definition}
\medskip

Let $i\rightarrow a_i\in\RR^+, \ i \in I$ be a function on $I$. Then, define $\sum\limits_{i\in I}a_i $ as the following element of $\RR^+\cup\{+\infty\}$:
\begin{eqnarray*}
\sum_{i\in I}a_i=\sup\bigg\{\sum_{i\in J}a_i: J\ \text{is a non-empty finite subset of }\ I\bigg\}.
\end{eqnarray*}
It is clear that, if $\sum\limits_{i\in I}a_i<+\infty $,
then the set of all $i\in I$ such that $a_i\neq0$ is at
most countable. In view of this, given a quaternionic Hilbert space $\H$ and a map $ i\rightarrow u_i\in\H$, $i\in I$,  the series $\sum\limits_{i\in I}u_i $ is said to be converges absolutely if $\sum\limits_{i\in I}\|u_i\|<+\infty$. If this happens, then only a finite or countable number of $u_i$ is non-zero and the series  $\sum\limits_{i\in I}u_i $ converges to a unique element of $\H$, independently from the ordering of  $u_i's$.

\begin{theorem}[\cite{GMP}]\label{2.6t}
	Let $\H$ be a quaternionic Hilbert space and let $N$ be a subset of $\H$
	such that, for $z, z'\in N$ such that $\langle z|z'\rangle=0$ if $z\neq z'$ and $\langle z|z\rangle=1$. Then the following conditions are equivalent:
	\begin{enumerate}[label=(\alph*)]
		\item For every $u,v\in\H$, the series $\sum_{z\in N}\langle u|z\rangle\langle z|v\rangle$ converges absolutely and
		\begin{eqnarray*}
		\langle u|v\rangle=\sum_{z\in N}\langle u|z\rangle\langle z|v\rangle.
			\end{eqnarray*}
			\item For every $u\in\H$, $\|u\|^2=\sum_{z\in N}|\langle z|u\rangle|^2$.
			\item $N^\bot={0}$.
			\item $\langle N\rangle$ is dense in $H$.
	\end{enumerate}
\end{theorem}

\begin{definition}[\cite{GMP}]
	Every quaternionic Hilbert space $\H$ admits a subset $N$, called \textit{Hilbert basis or orthonormal basis} of $\H$, such that, for $z, z'\in N$, $\langle z|z'\rangle=0$ if $z\neq z'$ and $\langle z|z\rangle=1$ and satisfies all the conditions of Theorem \ref{2.6t}.
	\end{definition}
Further, if there are two such sets, then they have the same cardinality. Furthermore, if $N$ is a Hilbert basis of $\H$, then for every $u\in\H$ can be uniquely expressed as
\begin{eqnarray*}
u=\sum_{z\in N}z\langle z|u\rangle
\end{eqnarray*}
where the series $\sum\limits_{z\in N}z\langle z|u\rangle$ converges absolutely in $\H$.
	 
\begin{definition}[\cite{adler}]
	Let $\H$ be a right quaternionic Hilbert space and $T$ be an operator on $\H$. Then $T$ is said to be
	\begin{itemize}[leftmargin=.25in]
		\item \emph{right $\HH$-linear} if
		$T( v_1\alpha + v_2\beta) = T(v_1)\alpha +  T(v_2)\beta, \ \text{for all} \ v_1, v_2 \in \H \ \text{and} \ \alpha, \beta \in \HH.$
		\item \emph{bounded} if there
		exist $K\ge 0$ such that
		$\pp{T(v)} \le K \pp{v}, \ \text{for all} \ v\in \H.$
	\end{itemize}
\end{definition}

\begin{definition}[\cite{adler}]
	Let $\H$ be a right quaternionic Hilbert space and $T$ be an operator on $\H$. Then the \textit{adjoint operator} $T^*$ of $T$ is defined by
	\begin{eqnarray*}
		\la v|Tu\ra = \la T^*v|u\ra, \ \text{for all} \ u, v \in\H.
	\end{eqnarray*}
	Further, $T$ is said to be \emph{self-adjoint} if $T=T^*$.
\end{definition}

\begin{theorem}[\cite{adler}]
	Let $\H$ be a right quaternionic Hilbert space and $S$ \& $T$ be two bounded right linear operators on $\H$. Then
	\begin{enumerate}[label=(\alph*)]
		\item $T+S$ and $TS\in\mathfrak{B}(\H)$. Moreover:
		\begin{eqnarray*}
		\|T+S\|\leq\|T\|+\|S\|\ \text{and}\ \|TS\|\leq\|T\|\|S\|
		\end{eqnarray*}
		\item $\la Tv|u\ra = \la v| T^*u \ra$.
		\item $(T+S)^*=T^*+S^*$.
		\item $(TS)^*=S^*T^*$.
		\item $(T^*)^*=T.$
		\item $I^*=I$, where $I$ is an identity operator on $\H$.
		\item If $T$ is an invertible operator then $(T^{-1})^*=(T^*)^{-1}$.
	\end{enumerate}
\end{theorem}

\begin{theorem}[\cite{GMP}]\label{2.11t}
	 Let $\H$ be a right quaternionic Hilbert space and let $T\in\mathfrak{B}(\H)$ be an operator. If $T\geq 0$,
	then there exists a unique operator in $\mathfrak{B}(\H)$, indicated by $\sqrt{T}$, such that $\sqrt{T}\geq 0$ 	and $\sqrt{T}\sqrt{T}=T$ . Furthermore, it turns out that $\sqrt{T}$ commutes with every operator
	which commutes with $T$ and if $T$ is invertible and self-adjoint, then $\sqrt{T}$ is also invertible and self-adjoint.	
\end{theorem}

\begin{theorem}[\cite{colombo}]
	Let $\H$ be a right quaternionic Hilbert space and  $T\in\mathfrak{B}(\H)$ and let $s\in\HH$ be such that $\|T\|<|s|$. Then the operator
\begin{eqnarray*}
	\sum_{n\geq 0}(s^{-1}T)^n{s^{-1}I}
\end{eqnarray*}
is the right and left algebraic inverse of $(sI-T)$ and the series converges in the operator norm.
\end{theorem}

\medskip

\section{Frames in Quaternionic Hilbert space}
\setcounter{equation}{0}
\fontsize{12}{14}

We begin this section with the following definition of frames in right quaternionic Hilbert spaces $\H$.
\begin{definition}Let $\H$ be a right quaternionic Hilbert space and $\{u_i\}_{i\in I}$ be a sequence in $\H$. Then $\cc{u_i}_{i\in I}$ is said to be a \textit{frame} for $V_R(\HH)$, if there exist two finite constants with $0<A\le B$  such that
	\begin{eqnarray}\label{3.1}
	A\|u\|^2\leq\sum_{i\in I}|\langle u_i|u\rangle|^2\leq B\|u\|^2, \ \text{for all}\ u\in V_R(\HH).
	\end{eqnarray}
	The positive constants $A$ and $B$, respectively, are called lower frame and upper frame bounds for the frame $\{u_i\}_{i\in I}$. The inequality (\ref{3.1}) is called frame inequality for the frame $\{u_i\}_{i\in I}$. A sequence $\{u_i\}_{i\in I}$ is called \textit{Bessel sequence} for right quaternionic Hilbert space $\H$ with bound $B$, if $\{u_i\}_{i\in I}$ satisfies the right hand side of the inequality (\ref{3.1}).
	A frame $\{u_i\}_{i\in I}$ for right quaternionic Hilbert space $V_R(\HH)$ is said to be
	\begin{itemize}[leftmargin=.25in]
		\item \textit{tight}, if it is possible to choose $A$ and $B$ satisfying inequality (3.1) with $A=B$.
		\item \textit{Parseval frame}, if it is tight with $A=B=1$.
		\item \textit{exact}, if it ceases to be a frame whenever anyone of its element is removed.
			\end{itemize}
	\end{definition}
Regarding the existence of frames in right quaternionic Hilbert space $\H$, we have the following examples:

\begin{example}
	Let $N$ be a Hilbert basis for right quaternionic Hilbert space  $\H$ such that, for each $z_i, z_k\in N$, $i, k\in I$, we have
\begin{eqnarray*}
	\langle z_i|z_k\rangle=
	\begin{cases}
		0,\ \text{for}\ i\neq k\\
		1,\ \text{for}\ i=k.
	\end{cases}
\end{eqnarray*}
\begin{enumerate}
	\item {Tight and Non-Exact}. Let $\{u_i\}_{i\in I}$ be a sequence in $\H$ defined as 
	\begin{eqnarray*}
	\qquad	u_i=u_{i+1}=z_i,\ i\in I.
	\end{eqnarray*}
Then $\{u_i\}_{i\in I}$ is tight and non-exact frame for $\H$ with bound $A=2$. 
Indeed, we have
\begin{eqnarray*}
\qquad\sum_{i\in I}|\langle u_i|u\rangle|^2&=&2\sum_{i\in I}|\langle z_i|u\rangle|^2\\&=&2\|u\|^2,\ \text{for all}\ u\in\H.
\end{eqnarray*}
\item {Non-Tight and Non-exact}. Let $\{u_i\}_{i\in I}$ be a sequence in $\H$ defined as 
\begin{eqnarray*}
	\qquad\begin{cases}
		u_1 = z_1\\
		u_i= z_{i-1},\ i\geq 2, i\in I.		
	\end{cases}
\end{eqnarray*}
Then $\{u_i\}_{i\in I}$ is non-tight and non-exact frame for $\H$. 
Indeed, we have
\begin{eqnarray*}
	\qquad\|u\|^2\leq\sum_{i\in I}|\langle u_i|u\rangle|^2\leq2\|u\|^2,\ \text{for all}\ u\in\H.
\end{eqnarray*}
\item {Parseval}. Let $\{u_i\}_{i\in I}$ be a sequence in $\H$ defined as
	\begin{eqnarray*}
	\qquad\begin{cases}
		u_1 = z_1\\
		u_{i_k}=u_{i_{k}+1}=u_{i_{k}+2}=\hdots =u_{i_{k+1}-1}=\dfrac{z_i}{\sqrt{i}},
	\end{cases}
\end{eqnarray*} 
$\text{where}\ i_k=i_{k-1}+(k-1),\ k\in\NN,\ i_0=1$.\\
Then $\{u_i\}_{i\in I}$ is a Parseval frame for $\H$. Indeed, we have
\begin{eqnarray*}
	\sum_{i\in I}|\langle u_i|u\rangle|^2=\sum_{i\in I}i\bigg|\bigg\langle \dfrac{z_i}{\sqrt{i}}\bigg|u\bigg\rangle\bigg|^2=\|u\|^2,\ \text{for all}\ u\in\H.
\end{eqnarray*}
\item {Exact}. Let $\{z_i\}_{i\in I}$ be a Hilbert basis of $\H$. Then $\{z_i\}_{i\in I}$ is an exact frame for $\H$.
	\end{enumerate}	
\end{example}

Next, we show  that for a sequence $\cc{u_i}_{i\in I}$ in right quaternionic Hilbert space $\H$ being a Bessel sequence is a sufficient condition for the series $\sum\limits_{i\in I} u_i q_i  $, $\{q_i\}_{i\in I}\subset \E$  to converge unconditionally.

\begin{theorem}\label{t3.3} Let $\H$ be a right quaternionic Hilbert space and  $\{u_i\}_{i\in I}$ be a Bessel sequence for $V_R(\HH)$ with  Bessel bound $B$. Then, for every sequence $\{q_i\}_{i\in I}\in\ell_2(\HH)$, the series $\sum\limits_{i\in I}u_iq_i$ converges unconditionally.
\end{theorem}

\proof For $i,j\in I,\ i>j$. Then, we have
\begin{eqnarray*}
	\bigg\|\sum_{k=1}^iu_kq_k-\sum_{k=1}^ju_kq_k\bigg\|&=&\bigg\|\sum_{k=j+1}^iu_kq_k\bigg\|\\
	&=&\sup_{\|v\|=1}\bigg|\bigg\langle\sum_{k=j+1}^iu_kq_k|v\bigg\rangle\bigg|\\
	&=&\sup_{\|v\|=1}\sum_{k=j+1}^i\bigg|\overline{q_k}\langle u_k|v\rangle\bigg|\\
	&\leq&\bigg(\sum_{k=j+1}^i|{q_k}|^2\bigg)^{1/2}\sup_{\|g\|=1}\bigg(\sum_{k=j+1}^i|\langle u_k|v\rangle|^2\bigg)^{1/2}\\
	&\leq&\sqrt{B}\bigg(\sum_{k=j+1}^i|{q_k}|^2\bigg)^{1/2}.
\end{eqnarray*}
Since $\{q_i\}_{i\in I}\in\ell_2(\HH)$ and $\left\{\sum\limits_{k=1}^i|q_k|^2\right\}_{i\in I}$ is a Cauchy sequence in $\RR$. Therefore $\left\{\sum\limits_{k=1}^iu_kq_k\right\}_{i\in I}$ is Cauchy sequence in $\H$. Hence $\left\{\sum\limits_{i\in I}u_iq_i\right\}$ unconditionally convergent in $\H$.
\endproof
\medskip

If view of Theorem \ref{t3.3}, if $\{u_i\}_{i\in I}$ is a Bessel sequence for $V_R(\HH)$. Then, the \textit{(right) synthesis operator} for $\{u_i\}_{i\in I}$ is a right linear operator $T:\ell_2(\HH)\to V_R(\HH)$ defined by
\begin{eqnarray*}T(\{q_i\}_{i\in I})=\sum_{i\in I}u_iq_i,\ \ \{q_i\}\in\ell_2(\HH).
\end{eqnarray*}
The adjoint operator $T^*$ of right synthesis operator $T$ is called the \textit{(right) analysis operator}. Further,  the  analysis operator $T^*:V_R(\HH)\to \ell_2(\HH)$  is given by
	\begin{eqnarray*}
		T^*(u)=\{\langle u_i|u\rangle\}_{i\in I},\ u\in V_R(\HH).
			\end{eqnarray*}
Infact, for $u\in V_R(\HH)$ and $\{q_i\}_{i\in I}\in\ell_2(\HH)$, we have
		\begin{eqnarray*}
			\langle T^*(u)|\{q_i\}_{i\in I}\rangle&=&\langle u|T(\{q_i\}_{i\in I})\rangle\\
			&=&\bigg\langle u\bigg|\sum_{i\in I}u_iq_i\bigg\rangle\\
			&=&\sum_{i\in I}\langle u|u_i\rangle q_i\\
			&=&\bigg\langle \{\langle u_i|u\rangle\}_{i\in I}, \{q_i\}_{i\in I}\bigg\rangle.
		\end{eqnarray*}
	Thus
	\begin{eqnarray*}
		T^*(u)=\{\langle u_i|u\rangle\}_{i\in I},\ u\in V_R(\HH).
	\end{eqnarray*}

Next, we give a characterization for a Bessel sequence in a right quaternionic Hilbert space.

\begin{theorem}\label{3.4t} Let $\H$ be a right quaternionic Hilbert space and  $\{u_i\}_{i\in I}$ be a sequence in $V_R(\HH)$. Then, $\{u_i\}_{i\in I}$ is a Bessel sequence for $V_R(\HH)$ with bound $B$ if and only if the right linear operator $T: \E \to \H$ defined by
\begin{eqnarray*}
T\left(\{q_i\}_{i\in I}\right)= \sum_{i\in I}u_iq_i,\ \{q_i\}_{i\in I}\in\ell_2(\HH)
\end{eqnarray*}
is a well defined and  bounded operator with $\|T\|\leq\sqrt{B}$.
\end{theorem}
\proof Let $\{u_i\}_{i\in I}$ be a Bessel sequence for right quaternionic Hilbert space $V_R(\HH)$. Then, by Theorem \ref{t3.3}, the $T$ is a well defined and  bounded operator with $\|T\|\leq \sqrt{B}$.

Conversely, let $T$ be a well defined and   bounded right linear operator with $\|T\|\leq\sqrt{B}$. Then, the adjoint of a bounded right linear operator $T$ is itself bounded and $\|T\|=\|T^*\|$. Since, for  $u\in V_R(\HH)$, we have
\begin{eqnarray*}
	\sum_{i\in I}|\langle u_i|u\rangle|^2&=&\|T^*(u)\|^2\\
	&\leq&\|T^*\|^2\|u\|^2\\
	&=&\|T\|^2\|u\|^2,
\end{eqnarray*} 
it follows that $\{u_i\}_{i\in I}$ is a Bessel sequence for right quaternionic Hilbert space $V_R(\HH)$ with bound $B$.
\endproof
\medskip

	Let $\H$ be a right quaternionic Hilbert space and $\{u_i\}_{i\in I}$ be a frame for $V_R(\HH)$. Then, the \textit{(right) frame operator} $S:V_R(\HH)\rightarrow V_R(\HH)$ for frame $\{u_i\}_{i\in I}$ is a right linear operator given by
	\begin{eqnarray*}
	S(u)&=&TT^*(u)\\
	&=&T(\{\langle u_i|u\rangle\}_{i\in I})\\
	&=&\sum_{i\in I}u_i\langle u_i|u\rangle,\ u\in \H.
	\end{eqnarray*}

\medskip

In next result, we discuss some properties of the  frame operator for a frame in right quaternionic Hilbert space.

\begin{theorem}
	Let $\H$ be a right quaternionic Hilbert space and  $\{u_i\}_{i\in I}$ be a frame for $V_R(\HH)$ with lower and upper frame bounds $A$ and $B$, respectively and  frame operator $S$. Then $S$ is positive, bounded, invertible and self adjoint right linear operator on $V_R(\HH)$. 
\end{theorem}
\proof For any $u\in V_R(\HH)$, we have
\begin{eqnarray*}
	\langle Su|u\rangle&=&\bigg\langle \sum_{i\in I}u_i\langle u_i|u\rangle|, u\bigg\rangle\\
	&=&\sum_{i\in I}\overline{\langle u_i|u\rangle}\langle u_i|u\rangle\\
	&=&\sum_{i\in I}|\langle u_i|u\rangle|^2.
\end{eqnarray*}
This gives
\begin{eqnarray*}
	A\|u\|^2\leq\langle Su|u\rangle\leq B\|u\|^2,\ u\in V_R(\HH).
\end{eqnarray*}
Thus 
\begin{eqnarray}\label{3.2}
AI\leq S\leq BI.
\end{eqnarray}
Hence $S$ is positive and bounded  right linear operator on $V_R(\HH)$. Also, $0\leq I-B^{-1}S\leq \dfrac{B-A}{B}I$ and consequently
\begin{eqnarray*}
	\|I-B^{-1}S\|&=&\sup_{\|v\|=1}\bigg|\bigg\langle (I-B^{-1}S)v\bigg|v\bigg\rangle\bigg|\\ &\leq& \dfrac{B-A}{B}\\&<&1.
\end{eqnarray*}
Then, $S$ is invertible. Further for any $u, v\in V_R(\HH)$, we have
\begin{eqnarray*}
\langle Su|v\rangle &=& \bigg\langle\sum_{i\in I}u_i\langle u_i|u\rangle\bigg| v\bigg\rangle\\
&=&\sum_{i\in I}\overline{\langle u_i|u\rangle}\langle u_i|v\rangle\\
&=&\sum_{i\in I}\langle u|u_i\rangle\langle u_i|v\rangle\\
\qquad\quad&=&\bigg\langle u\bigg|\sum_{i\in I}u_i\langle u_i|v\rangle\bigg\rangle\\
&=&\langle u|Sv\rangle.
\end{eqnarray*}
Thus $S$ is also a self adjoint right linear operator on $\H$.\endproof

\begin{corollary}[The Reconstruction Formula]
		Let $\H$ be a right quaternionic Hilbert space and  $\{u_i\}_{i\in I}$ be a frame for $V_R(\HH)$ with  frame operator $S$. Then, for every $u\in \H$ can be expressed as
		\begin{eqnarray*}
		u=\sum_{i\in I}S^{-1}u_i\langle u_i|u\rangle.
	\end{eqnarray*}
\end{corollary}

\proof As $S$ is invertible. Therefore, for $u\in \H$, we have
\begin{eqnarray*}
	u&=&S^{-1}S(u)\\
	&=&\sum_{i\in I}S^{-1}u_i\langle u_i|u\rangle.\qquad\qquad\qquad\qquad\qquad\qquad\qquad\qquad\qquad\qquad\qquad\qquad \Box
\end{eqnarray*}

In the next result, we construct a frame with the help a give frame in a right quaternionic Hilbert space.
\begin{theorem}\label{3.7t}
		Let $\H$ be a right quaternionic Hilbert space and  $\{u_i\}_{i\in I}$ be a frame for $V_R(\HH)$ with lower and upper frame bounds $A$ and $B$, respectively and  frame operator $S$.  Then $\{S^{-1}u_i\}_{i\in I}$ is also a frame for $\H$ with bounds $B^{-1}$ and ${A^{-1}}$ and right frame operator $S^{-1}$. 
\end{theorem}
\proof For   $u\in\H$, we have
\begin{eqnarray*}
	\sum_{i\in I}|\langle S^{-1}u_i|u\rangle|^2&=&\sum_{i\in I}|\langle u_i|S^{-1}u\rangle|^2\\
	&\leq&B\|S^{-1}u\|^2\\
	&\leq&B\|S^{-1}\|^2\|u\|^2.
\end{eqnarray*}
Thus $\{S^{-1}u_i\}_{i\in I}$ is a Bessel sequence for $\H$. It follows that the right frame operator for $\{S^{-1}u_i\}_{i\in I}$ is well defined. Therefore, we have
\begin{eqnarray}\label{3.3}
\sum_{i\in I}S^{-1}u_i\langle S^{-1}u_i|u\rangle&=&S^{-1}\bigg(\sum_{i\in I}u_i\langle S^{-1}u_i|u\rangle\bigg) \nonumber\\
&=&S^{-1}S(S^{-1}u)\nonumber\\&=&S^{-1}u,\ u\in\H.
\end{eqnarray}
Thus, the right frame operator for $\{S^{-1}u_i\}_{i\in I}$ is $S^{-1}$. The operator $S^{-1}$ commutes with both $S$  and $I:\H\rightarrow \H$ (an identity operator on $\H$). Therefore, multiplying  the  inequality (\ref{3.2}) with $S^{-1}$, we have
\begin{eqnarray*}
&&	B^{-1}I\leq S^{-1}\leq A^{-1}I\\
\text{i.e.}&&	B^{-1}\|u\|^2\leq\langle S^{-1}u|u\rangle\leq A^{-1}\|u\|^2,\ u\in\H.
\end{eqnarray*}
By using (\ref{3.3}), we get
\begin{eqnarray*}
B^{-1}\|u\|^2\leq\sum_{i\in I}|\langle S^{-1}u_i|u\rangle|^2\leq A^{-1}\|u\|^2,\ u\in\H.
\end{eqnarray*}
Therefore $\{S^{-1}u_i\}_{i\in I}$ is frame for $\H$ with bounds $\dfrac{1}{B}$ and $\dfrac{1}{A}$ and right frame operator $S^{-1}$.
\endproof

In view of Theorem \ref{3.7t}, we have a following definition:

\begin{definition}
	Let $\H$ be a right quaternionic Hilbert space and $\{u_i\}_{i\in I}$ be a frame for  $\H$ with  frame operator $S$. Then, the frame $\{S^{-1}u_i\}_{i\in I}$ is called the \textit{canonical dual} of the frame $\{u_i\}_{i\in I}$.
\end{definition}

Next, we construct a Parseval frame with the help of a give frame in a right quaternionic Hilbert space.
\begin{theorem}
		Let $\H$ be a right quaternionic Hilbert space and  $\{u_i\}_{i\in I}$ be a frame for $V_R(\HH)$ with   frame operator $S$. Then $\{S^{-1/2}u_i\}_{i\in I}$ is  a Parseval frame for $\H$. 
\end{theorem}
\proof By Theorem \ref{2.11t}, for any $u\in\H$, we have
\begin{eqnarray}\label{3.4}
	u&=&S^{-1/2}SS^{-1/2}u\nonumber\\
	&=&S^{-1/2}\sum_{i\in I}u_i\langle u_i|S^{-1/2}u\rangle\nonumber\\
	&=&\sum_{i\in I}S^{-1/2}u_i\langle u_i|S^{-1/2}u\rangle.
\end{eqnarray}
Therefore, for any $u\in\H$
\begin{eqnarray*}
	\|u\|^2&=&\langle u|u\rangle\\
	&=&\bigg\langle\sum_{i\in I}S^{-1/2}u_i\langle u_i|S^{-1/2}u\rangle\bigg|u\bigg\rangle\ \ \ \ \ \left(\text{by using (\ref{3.4})}\right).\\
	&=&\sum_{i\in I}\overline{\langle u_i|S^{-1/2}u\rangle}\langle S^{-1/2}u_i|u\rangle\\
	&=&\sum_{i\in I}|\langle S^{-1/2}u_i|u\rangle|^2.
\end{eqnarray*}
Thus $\{S^{-1/2}u_i\}_{i\in I}$ is a Parseval frame for $\H$.
\endproof

Next, we give a characterization of Parseval frames $\{u_i\}_{i\in I}$ for right quaternionic Hilbert space $\H$.

\begin{theorem}
	Let $\H$ be a right quaternionic Hilbert space and $\{u_i\}_{i\in I}$ be a frame for  $\H$ with  frame operator $S$. Then $\{u_i\}_{i\in I}$ is a Parseval frame for $\H$ if and only if $S$ is an identity operator on $\H$.
\end{theorem}
\proof Let $\{u_i\}_{i\in I}$ be a Parseval frame for $\H$. Then, for all $u\in \H$, we have
\begin{eqnarray*}
	\sum_{i\in I}|\langle u_i|u\rangle|^2&=&\|u\|^2
\end{eqnarray*}
which implies
\begin{eqnarray*}
	\langle Su|u\rangle=\langle u|u\rangle.
\end{eqnarray*}
Thus, $S$ is an identity operator on $\H$.

Conversely, let $S$ is an identity operator on $\H$. Thus, for  $u\in\H$, we have
\begin{eqnarray}\label{3.5}
u&=&S(u)\nonumber\\&=&\sum_{i\in I}u_i\langle u_i|u\rangle.
\end{eqnarray}
So, we have
\begin{eqnarray*}
	\|u\|^2&=&\langle u|u\rangle\\
	&=&\bigg\langle\sum_{i\in I}u_i\langle u_i|u\rangle|u\bigg\rangle \\
	&=&\sum_{i\in I}|\langle u_i|u\rangle|^2.
\end{eqnarray*}
Hence $\{u_i\}_{i\in I}$ is a Parseval frame for right quaternionic Hilbert space $\H$.
\endproof

In the next theorem, we give a necessary condition for a Parseval frame $\{u_i\}_{i\in I}$ for right quaternionic Hilbert space $\H$.

\begin{theorem}
		Let $\H$ be a right quaternionic Hilbert space and  $\{q_i\}_{i\in I}$ be a semi-normalized sequence of quaternions in $\HH$ with bounds $a$ and $b$. If $\{u_iq_i\}_{i\in I}$ is a Parseval frame for $\H$, then $\{u_i\}_{i\in I}$ is a frame for $\H$ with bounds ${b^{-2}}$ and ${a^{-2}}$.
\end{theorem}
\proof Let $\{u_iq_i\}_{i\in I}$ be a Parseval frame for $\H$. Then, for any $u\in\H$, we have
\begin{eqnarray*}
\sum_{i\in I}|\langle u_iq_i|u\rangle|^2=\|u\|^2.
\end{eqnarray*}
This gives
\begin{eqnarray}
\sum_{i\in I}|q_i|^2|\langle u_i|u\rangle|^2=\|u\|^2.
\end{eqnarray}
Since $\{q_i\}_{i\in I}$ is a normalized sequence with bounds $a$ and $b$, we get
\begin{eqnarray*}
a^2\sum_{i\in I}|\langle u_i|u\rangle|^2\leq\sum_{i\in I}|q_i|^2|\langle u_i|u\rangle|^2\leq b^2\sum_{i\in I}|\langle u_i|u\rangle|^2,\ u\in\H.
\end{eqnarray*}
Therefore, by using (\ref{3.5}), we have
\begin{eqnarray*}
	\dfrac{1}{b^2}\|u\|^2\leq\sum_{i\in I}|\langle u_i|u\rangle|^2\leq\dfrac{1}{a^2},\ u\in\H.
\end{eqnarray*}
Hence $\{u_i\}_{i\in I}$ is a frame for $\H$ with bounds ${b^{-2}}$ and ${a^{-2}}$.
\endproof

Finally, in this section we give a characterization of frame for right quaternionic Hilbert space in terms of operators.

\begin{theorem}	Let $\H$ be a right quaternionic Hilbert space and $\{u_i\}_{i\in I}$ be a sequence in  $\H$. Then $\{u_i\}_{i\in I}$ is a frame for $\H$ if and only if the right linear operator $T:\E\to \H$ 
	\begin{eqnarray*}
	T(\{q_i\}_{i\in I})=\sum_{i\in I}u_iq_i, \  \ \cc{q_i}_{i\in I} \in \H
	\end{eqnarray*}
	is a well-defined and bounded mapping from $\ell_2(\HH)$ onto $\H$.
\end{theorem}
\proof Let $\{u_i\}_{i\in I}$ be a frame for $\H$. Then the  frame operator $S=TT^*$
for $\{u_i\}_{i\in I}$ is invertible on $\H$. So $T$ is an onto mapping. Also, by Theorem \ref{3.4t}, $T$ is well defined and bounded mapping from $\ell_2(\HH)$ to $\H$.

Conversely, let $T$ be a well defined, bounded right linear operator from $\ell_2(\HH)$ onto $\H$. Then, by Theorem \ref{3.4t}, $\{u_i\}_{i\in I}$ is a Bessel sequence for $\H$. Since $T$ is an onto, there exists an right linear operator $T^{\dag}:\H\rightarrow\ell_2(\HH)$ such that
\begin{eqnarray*}
u&=&TT^{\dag}u\\
&=&\sum_{i\in I}u_i(T^{\dag}u)_i,\ \ u\in\H,
\end{eqnarray*}
where $(T^{\dag}u)_i$ denotes the $i-$th coordinate of $T^{\dag}u$.
Therefore, for $u\in\H$, we have
\begin{eqnarray*}
\|u\|^4&=&|\langle u|u\rangle|^2\\
&=&\bigg|\bigg\langle\sum_{i\in I}u_i(T^{\dag}u)_i\bigg|u\bigg\rangle\bigg|^2\\
&\leq&\sum_{i\in I}|(T^{\dag}u)_i|^2\ \sum_{i\in I}|\langle u_i|u\rangle|^2\\
&\leq&\|T^{\dag}\|^2\|u\|^2\sum_{i\in I}|\langle u_i|u\rangle|^2.
\end{eqnarray*}
Hence $\{u_i\}_{i\in I}$ is a frame for $\H$ with bounds $\|T^{\dag}\|^{-2}$ and $\|T\|^2$. 
\endproof

\section{Stability of Frames in Right Quaternionic Hilbert Space}
\setcounter{equation}{0}
\fontsize{12}{14}

Christensen \cite{CH1} gave a version of Paley-Wiener Theorem for frames in Hilbert spaces. In this section, we give a similar version for frames in right quaternionic Hilbert spaces.

\begin{theorem}\label{4.1t}
	Let $\H$ be a right quaternionic Hilbert space and $\{u_i\}_{i\in I}$ be a frame for $\H$ with lower and upper frame bounds $A$ and $B$, respectively. Let $\{v_i\}_{i\in I}$ be a sequence  in $\H$ and assume that there exist $\lambda,\ \mu\geq 0$ such that $\bigg(\lambda+\dfrac{\mu}{\sqrt{A}}\bigg)<1$ and 
	\begin{eqnarray}\label{4.1}
	\bigg\|\sum_{i\in J}(u_i-v_i)q_i\bigg\|\leq\lambda\bigg\|\sum_{i\in J}u_iq_i\bigg\|+\mu\bigg(\sum_{i\in J}|q_i|^2\bigg)^{1/2}
	\end{eqnarray}
	for all finite quaternions $q_i\in\HH,\ i\in J\subseteq I$ with $|J|<+\infty$. Then $\{v_i\}_{i\in I}$ is a frame for $\H$ with bounds $A\bigg(1-\bigg(\lambda+\dfrac{\mu}{\sqrt{A}}\bigg)\bigg)^2$
	and $B\bigg(1+\bigg(\lambda+\dfrac{\mu}{\sqrt{B}}\bigg)\bigg)^2$.
\end{theorem}

\proof Let $J\subseteq I$ with $|J|<+\infty$. Then
\begin{eqnarray*}
	\bigg\|\sum_{i\in J}v_iq_i\bigg\|&\leq&\bigg\|\sum_{i\in J}(u_i-v_i)q_i\bigg\|+\bigg\|\sum_{i\in J}u_iq_i\bigg\|\\
	&\leq&(1+\lambda)\bigg\|\sum_{i\in J}u_iq_i\bigg\|+\mu\bigg(\sum_{i\in J}|q_i|^2\bigg)^{1/2}.
\end{eqnarray*}
Since
\begin{eqnarray*}
	\bigg\|\sum_{i\in J}u_iq_i\bigg\|&=&\sup_{\|v\|=1}\bigg|\bigg\langle\sum_{i\in J}u_iq_i\bigg|v\bigg\rangle\bigg|\\
	&\leq&\bigg(\sum_{i\in J}|q_i|^2\bigg)^{1/2}\sup_{\|v\|=1}\bigg(\sum_{i\in J}|\langle u_i|v\rangle|^2\bigg)^{1/2}\\
	&\leq&\sqrt{B}\bigg(\sum_{i\in J}|q_i|^2\bigg)^{1/2}.
\end{eqnarray*}
This gives
\begin{eqnarray}\label{4.2}
\bigg\|\sum_{i\in J}v_iq_i\bigg\|\leq\bigg({(1+\lambda)}\sqrt{B}+{\mu}\bigg)\bigg(\sum_{i\in J}|q_i|^2\bigg)^{1/2}.
\end{eqnarray}
Let $T:\ell_2(\HH)\rightarrow\H$ be a right linear operator defined by 
\begin{eqnarray*}
	T(\{q_i\}_{i\in I})=\sum_{i\in J}v_iq_i,\ \{q_i\}_{i\in I}\in\ell_2(\HH).
\end{eqnarray*}
Then, by (\ref{4.2}), $T$ is well defined and bounded with $\|T\|\leq\bigg({(1+\lambda)}\sqrt{B}+{\mu}\bigg)$. Therefore, by Theorem \ref{3.4t}, $\{v_i\}_{i\in I}$ is a Bessel sequence for $\H$ with the upper bound $B\bigg(1+\bigg(\lambda+\dfrac{\mu}{\sqrt{B}}\bigg)\bigg)^2$.\\
Let $U$ be the right synthesis operator for the frame $\{u_i\}_{i\in I}$.
Define an operator $W:\H\rightarrow\ell_2(\HH)$ by
\begin{eqnarray*}
	W(u)=U^*(UU^*)^{-1}u,\ u\in\H.
\end{eqnarray*}
Since $\{(UU^*)^{-1}u_i\}_{i\in I}$ is also a frame for $\H$ with bounds $\dfrac{1}{B}$ and $\dfrac{1}{A}$, then
\begin{eqnarray*}
	\|W(u)\|^2&=&\sum_{i\in I}|\langle S^{-1}u_i|u\rangle|^2\\
	&\leq&\dfrac{1}{A}\|u\|^2,\ u\in\H.
\end{eqnarray*}
Using (\ref{4.1}) with $\{q_i\}_{i\in I}=Wu,\ u\in\H$, we get
\begin{eqnarray*}
	\|u-TWu\|&\leq&\lambda\|u\|+\mu\|Wu\|\\
	&\leq&\bigg(\lambda+\dfrac{\mu}{\sqrt{A}}\bigg)\|u\|,\ u\in\H.
\end{eqnarray*}
This gives $\|TW\|\leq 1+\lambda+\dfrac{\mu}{\sqrt{A}}$. Since $\bigg(\lambda+\dfrac{\mu}{\sqrt{A}}\bigg)<1$, then the operator $TW$ is invertible and $\|(TW)^{-1}\|\leq\dfrac{1}{1-\bigg(\lambda+\dfrac{\mu}{\sqrt{A}}\bigg)}$.

\noindent
Now, for any $u\in\H$, we have
\begin{eqnarray*}
	u&=&(TW)(TW)^{-1}u\\
	&=&\sum_{i\in I}v_i\langle (UU^*)^{-1}u_i|(TW)^{-1}u\rangle.
\end{eqnarray*}
Therefore, for each $u\in\H$, we have
\begin{eqnarray*}
	\|u\|^4&=&\bigg|\bigg\langle \sum_{i\in I}v_i\langle (UU^*)^{-1}u_i|(TW)^{-1}u\rangle\bigg|u\bigg\rangle\bigg|^2\\
	&\leq&\sum_{i\in I}\bigg|\langle (UU^*)^{-1}u_i|(TW)^{-1}u\rangle\bigg|^2\bigg(\sum_{i\in I}|\langle v_i|u\rangle|^2\bigg)\\
	&\leq&\dfrac{1}{A}\bigg[ \dfrac{1}{1-(\lambda+\frac{\mu}{\sqrt A})}\bigg]^2\|u\|^2\bigg(\sum_{i\in I}|\langle v_i|u\rangle|^2\bigg).
\end{eqnarray*}
This gives
\begin{eqnarray*}
	\sum_{i\in I}|\langle v_i|u\rangle|^2\geq A\bigg(1-\bigg(\lambda+\dfrac{\mu}{\sqrt{A}}\bigg)\bigg)^2 \|u\|^2,\ u\in\H.
\end{eqnarray*}
Hence $\cc{v_i}_{i\in I}$ is a frame for $\H$ with desire bounds.
\endproof

\begin{corollary}
Let $\H$ be a right quaternionic Hilbert space and $\{u_i\}_{i\in I}$ be a frame for $\H$ with lower and upper frame bounds $A$ and $B$, respectively. Let $\{v_i\}_{i\in I}$ be a sequence in $\H$ and assume that there exist $0<R<A$ such that
	\begin{eqnarray}\label{4.3}
	\|\sum_{i\in I}(u_i-v_i)q_i\|\leq\sqrt{R}\bigg(\sum_{i\in I}|q_i|^2\bigg)^{1/2},\ \text{for all}\ \{q_i\}_{i\in I}\in\ell_2(\HH).
	\end{eqnarray}
	Then $\{v_i\}_{i\in I}$ is also a frame for $\H$ with bounds $(\sqrt{A}-\sqrt{R})^2$ and $(\sqrt{B}+\sqrt{R})^2$.
\end{corollary}
\proof Take $\lambda=0$ and $\mu=\sqrt{R}$ in Theorem \ref{4.1t}.
\endproof

\begin{remark}
	A condition $\bigg(\lambda+\dfrac{\mu}{\sqrt{A}}\bigg)<1$ can not be dropped in the Theorem \ref{4.1t}. In this, regard, we give the following example:
\end{remark}

\begin{example}
	Let $\H$ be a right quaternionic Hilbert space, $\{z_i\}_{i\in \NN}$ be a Hilbert basis of $\H$, $\{p_i\}_{i\in\NN}$ be a sequence in $\HH$ and $\{v_i\}_{i\in\NN}$ be a sequence in $\H$ such that
	\begin{eqnarray*}
		v_i=z_i+z_{i+1}p_i,\ i\in\NN.
	\end{eqnarray*}
	Then, for $J\subseteq\NN$ with $|J|<\infty$ and $\{q_i\}_{i\in\NN}\in\ell_2(\HH)$, we have
	\begin{eqnarray*}
		\bigg\|\sum_{i\in J}(v_i-z_i)q_i\bigg\|&=&\bigg\|\sum_{i\in J}z_{i+1}p_iq_i\bigg\|\\
		&\leq&\sup_{i}\ |p_i|\bigg(\sum_{i\in I}|q_i|^2\bigg)^{1/2}.
	\end{eqnarray*}
	Thus, if $p=\sup\limits_{i} \ |p_i|<1$, then by Theorem \ref{4.1t}, $\{u_i\}_{i\in\NN}$ is a frame for $\H$ with bounds $(1-p)^2$ and $(1+p)^2$.	On taking $p_i=1$, for all $i\in\NN$, we get
	\begin{eqnarray}\label{4.4}
	v_i=z_i+z_{i+1},\ i\in\NN.
	\end{eqnarray}
	Then, for any sequence $\{q_i\}_{i\in\NN}\in\ell_2(\HH)$ and for $J\subseteq\NN$ with $|J|<\infty$ 
	\begin{eqnarray*}
		\bigg\|\sum_{i\in J}(v_i-z_i)q_i\bigg\|&=&\bigg\|\sum_{i\in J}z_{i+1}q_i\bigg\|\\
		&\leq&\bigg(\sum_{i\in I}|q_i|^2\bigg)^{1/2}.
	\end{eqnarray*}
	Thus, the condition (\ref{4.1}) is satisfied with either $(\lambda, \mu)=(1,0)\ \text{or}\ (0,1)$. 
	So, $\bigg(\lambda+\dfrac{\mu}{\sqrt{A}}\bigg)\nless1$ and $\{v_i\}_{i\in\NN}$ is not a frame for $\H$.
\end{example}

\end{document}